\documentclass[reqno]{amsart}
\usepackage{amsmath}
\usepackage{amssymb}
\usepackage{amsfonts}
\usepackage{amstext}
\usepackage[cp866]{inputenc}
\usepackage[english]{babel}
\theoremstyle{plain}\newtheorem{lem}[]{Lemma}
\theoremstyle{plain}
\theoremstyle{plain}
\theoremstyle{plain}\newtheorem{cor}[]{\bf{Corollary}}

\title
[The distribution of random evolution...] {The distribution of
random evolution in Erlang semi-Marov media}

\author{A. Pogorui}
\address{Shelushkova 10, 222, Zhitomir, Ukraine}
\curraddr{Shelushkova 10, 222, Zhitomir, Ukraine }
\email{pogor@zu.edu.ua}


\begin{document}

%    Information for second author
%\author{U. T. Hora}
%\address{Address of U. T. Hora}
%\curraddr{Current address of U. T. Hora}
%\email{email@ of U. T. Hora}
%\thanks{This is thanks of U. T. Hora}

%    Information for third author
%\author{T. H. Orau}
%\address{Address of T. H. Orau}
%\curraddr{Current address of T. H. Orau}
%\email{email@ of T. H. Orau}
%\thanks{This is thanks of T. H. Orau}

%    General info
\date{01/01/2009}
%\dedicatory{}

\begin{abstract}
In this paper we study a one-dimensional random motion by having a
general Erlang distribution for the sojourn times of the switching
process and we obtain solution of the four order hyperbolic PDE
for 2-Erlang case.
\end{abstract}

\maketitle

\section{Introduction}

%\subsection{This is a subsection}

In the paper [1] we studied a one-dimensional random motion with
the $m$-Erlang distribution between consequent epochs of velocity
alternations. Let $f(t,x)$ be the probability density function
(pdf) of a particle position at time $t$, provided that it exists.
We obtained the following higher order hyperbolic equations for
$f(t,x)$

\begin{equation}\label{1}
\left( \frac{\partial}{\partial t}-v\frac{\partial}{\partial
x}+\lambda \right)^{m}\left(\frac{\partial}{\partial
t}+v\frac{\partial}{\partial x}+\lambda
\right)^{m}f(t,x)-\lambda^{2m}f(t,x)=0,
\end{equation}
where $v>0$ is the velocity of the particle and $\lambda$ is the
parameter of the $m$-Erlang distribution. It is assumed that a
particle started at $x=0$ and hence, $f(0,x)=\delta(x)$.

The pdf $f(t,x)$ can be represented in the following form
$f(t,x)=f_c (t,x)+f_s (t,x)$, where $f_c (t,x)$ is the absolute
continuous part and $f_s (t,x)$ is the singular part w.r.t.
Lebesgue measure on the line.

\begin{lem}
The singular part $f_s (t,x)$ of the pdf $f(t,x)$ is of the
following form
\begin{equation}\label{2}
f_s (t,x/v)=\delta(t-x/v)e^{-\lambda t} \sum_{i=0}^{m-1}(\lambda
t)^i/i!,
\end{equation}
\end{lem}

\emph{Proof.} It is evident that for $t=x/v$ the pdf $f(t,x)$ has
the singularity given by Eq.(2). Let us show that for $t>|x/v|$
the pdf $f(t,x)$ has no singularity w.r.t. Lebesgue measure on
$\mathbb{R}$. Denote by $v_k$ the random event "$k$ velocity
alternations occurred". For $\Delta x=[x,x+\Delta]$, $\Delta>0$,
consider
\begin{equation}
P_{\bar{\nu}_0} (x(t)\in \Delta x)=\sum_{k\geq 1} P(x(t)\in \Delta
x,\nu_k ),\nonumber
\end{equation}
which is the probability of the event that at least one
alternation occurred and $x(t)\in \Delta x$. Let us show that for
each $t>0$ there exists a constant $C_t<\infty$ such that
\begin{equation}
\sup_{x} \frac{P_{\bar{\nu}_0 } (x(t)\in \Delta x)}{\Delta x}
<C_t. \nonumber
\end{equation}

Denote by $\theta_k, k\geq 1,$ time between $(k-1)$th and $k$th
velocity alternations. Recall that $\theta_k$,  $k \geq 1$ are
independent $m$-Erlang distributed random variables. It is easily
verified that

\begin{align}
P_{\bar{\nu}_0 } (x(t)\in \Delta x) = & \sum_{k\geq 1}
P\left(\sum_{i=1}^{k}(-1)^{i+1} \theta_i v +(-1)^k
\left(t-\sum_{i=1}^{k} \theta_i v \right) \in \Delta x,
\sum_{i=1}^k \theta_i <t
\right) \nonumber\\
&= \sum_{k\geq 1} P\left(\left(\sum_{i=1}^{k}(-1)^{i+1} \theta_i
-(-1)^k \sum_{i=1}^{k} \theta_i \right)v \in \Delta x-(-1)^{k}vt,
\sum_{i=1}^k \theta_i <t
\right) \nonumber \\
&=\sum_{l\geq 0}P\left(2v(\theta_1+ \theta_3+\ldots+\theta_{2l+1})
\in \Delta x-vt, \sum_{i=1}^{2l+1}
\theta_i <t \right) \nonumber \\
&= \sum_{l\geq 0}P\left(-2v(\theta_2+
\theta_4+\ldots+\theta_{2l+2}) \in \Delta x+vt, \sum_{i=1}^{2l+2}
\theta_i <t \right)\nonumber \\
&\leq \sup_{x}\sum_{l\geq 0}P\left(2v\sum_{i=1}^{l} \theta_{2i-1}
\in \Delta x, 2v\sum_{i=1}^{l}
\theta_{2i} <vt-x\right) \nonumber \\
&+ \sup_{x}\sum_{l\geq 0}P\left(-2v\sum_{i=1}^{l} \theta_{2i} \in
\Delta x, 2v\sum_{i=1}^{l} \theta_{2i+1} <vt+x\right).\nonumber
\end{align}

Since $|x|\leq vt$ and for every $m\geq 1$ the pdf $p_m
(x,\lambda)$ of the $m$-Erlang distribution with the parameter
$\lambda$ satisfies $p_m (x,\lambda)\leq \lambda$, we have

\begin{align}\label{3}
\sum_{l\geq 1} P(2v(\theta_1+\theta_3+\ldots+\theta_{2l-1})\in
\Delta x,2v(\theta_2+\theta_4+\ldots+\theta_{2l} )< vt-x)\nonumber \\
\leq \frac{\lambda\Delta}{2v} \sum_{l\geq
1}P(\theta_2+\theta_4+\ldots+\theta_{2l}<t)
\end{align}

Since $\theta_i$ is $m$-Erlang distributed we have for $2lm+1>t$
\begin{equation*}
P(\theta_2+\theta_4+\ldots+\theta_{2l}<t)\leq \left(e^{\lambda
t}-\sum_{i=0}^{2lm}\frac{(\lambda t)^{i}}{i!}\right) e^{-\lambda
t}\leq \frac{(\lambda t)^{2lm+1} e^{-\lambda
t}}{2lm!(2lm+1-\lambda t)}.
\end{equation*}

Therefore, taking into account (3), there exists a constant $A_t$
such that
\begin{equation*}
\sup_{x}\sum_{l\geq 1}P\left(2v\sum_{i=1}^{l} \theta_{2i-1} \in
\Delta x, 2v\sum_{i=1}^{l} \theta_{2i} <vt-x\right)\leq A_t
\Delta.
\end{equation*}

In much the same way, we can show that there exists a constant
$B_t$ such that

\begin{equation*}
\sup_{x}\sum_{l\geq 1}P\left(-2v\sum_{i=1}^{l} \theta_{2i} \in
\Delta x, 2v\sum_{i=1}^{2l-1} \theta_{2i} <vt+x\right)\leq B_t
\Delta.
\end{equation*}

Putting $C_t=A_t+B_t$, we conclude the proof.

\begin{cor}
The absolute continuous part $f_c (t,x)$ of the pdf $f(t,x)$
satisfies Eq.(1) for $t<|\frac{x}{v}|$.
\end{cor}

Now let us study the behavior of the continuous part $f_c (t,x)$
close to lines $t=\pm\frac{x}{v}$.

\begin{lem}
For $m\geq 2$, we have
\begin{align}
\lim_{\varepsilon \downarrow 0}
\frac{P\{0<t-x(t)<\varepsilon\}}{\varepsilon} =\frac{\lambda^m
t^{m-1}
e^{-\lambda t}}{2(m-1)!},\nonumber \\
\lim_{\varepsilon\rightarrow 0}
\frac{P\{t+x(t)<\varepsilon\}}{\varepsilon}=0.\nonumber
\end{align}
\end{lem}

\emph{Proof.} It is easily verified that
\begin{equation*}
P\left\{0<t-x(t)\leq \varepsilon
\right\}=P\left\{t-\frac{\varepsilon}{2}\leq \theta_1<t
\right\}+\int_0^tP\left\{\theta_3\geq t-u,\theta_2\leq
\frac{\varepsilon}{2}, \theta_1\in du \right\} +o(\varepsilon),
\end{equation*}
where $\theta_i, i=1,2,3$ are independent $m$-Erlang distributed
random variables with the parameter $\lambda$. Since $\int_0^t
P(\theta_3\geq t-u,\theta_2\leq \frac{\varepsilon}{2},\theta_1\in
du)=o(\varepsilon)$, passing to the limit, we get
\begin{align}
\lim_{\varepsilon \downarrow 0}\
\frac{P\left\{0<t-x\left(t\right)<\varepsilon
\right\}}{\varepsilon }\nonumber \\= & \lim_{\varepsilon
\downarrow 0}\ \frac{e^{-\lambda t}}{\varepsilon
}\left(\left(\sum^{m-1}_{i=0}{\frac{{\left(\lambda
t\right)}^i}{i!}}\right)-e^{\lambda \frac{\varepsilon
}{2}}\left(\sum^{m-1}_{i=0}{\frac{{\left(\lambda
\left(t-\frac{\varepsilon
}{2}\right)\right)}^i}{i!}}\right)\right)\nonumber \\  &
=\frac{{\lambda }^mt^{m-1}e^{-\lambda
t}}{2\left(m-1\right)!}\nonumber.
\end{align}

Similarly, $P\left\{t+x\left(t\right)\le \varepsilon
\right\}=P\left\{{\theta }_2\ge t-\frac{\varepsilon }{2},{\theta
}_1\le \frac{\varepsilon }{2}\right\}+o\left(\varepsilon \right)$
and as it easily seen that
\[\mathop{{\rm lim}}_{\varepsilon \downarrow 0}\ \frac{P\left\{t+x\left(t\right)<\varepsilon \right\}}{\varepsilon }=0.\]

%\begin{rem}
The case where $m=1$ will be considered below as an example.
%\end{rem}

%\begin{rem}
We will seek solutions of Eq.(1) among functions which continuous
part $f_c (t,x)$ satisfies the following conditions

\begin{equation}\label{4}
\lim_{x\uparrow t} f_c\left(t,x\right)\ =\lim_{\varepsilon
\downarrow 0}\ \frac{P\left\{0<t-x\left(t\right)<\varepsilon
\right\}}{\varepsilon },\\
\lim_{x\downarrow -t} f_c\left(t,x\right)\ \lim_{\varepsilon
\downarrow 0}\ \frac{P\left\{t+x\left(t\right)<\varepsilon
\right\}}{\varepsilon }.
\end{equation}
%\end{rem}

\noindent By applying the transformation
$f\left(t,x\right)=e^{\lambda t}g\left(t,x\right)$ and changing
the variable $y=\frac{x}{v}$, we reduce Eq.(1) to

\begin{equation}\label{5}
\left(\frac{{\partial }^2}{\partial t^2}-\frac{{\partial
}^2}{\partial y^2}\right)^m g_c\left(t,y\right)-{\lambda
}^{2m}g_c\left(t,y\right)=0,
\end{equation}
with the singular part
$g_s\left(t,y\right)=\left(\sum^{m-1}_{i=0}{\frac{{\left(\lambda
t\right)}^i}{i!}}\right)\delta \left(t-y\right).$

In the sequel we assume, without restricting the generality, that
$\lambda =1$. By introducing the function ${\mathbf
f}\left(t,y,z\right)=e^zg_c(t,y)$, we reduce Eq.(5) to the
following equation

\begin{equation}\label{6}
\left(\frac{{\partial }^2}{\partial t^2}-\frac{{\partial
}^2}{\partial y^2}\right)^m{\mathbf
f}\left(t,y,z\right)-\frac{{\partial }^{2m}}{\partial
z^{2m}}{\mathbf f}\left(t,y,z\right)=0.
\end{equation}

\noindent We will seek solutions of this equation by using theory
of differentiable functions on commutative algebras [2].

\section{Main results}

Let $A_0$ be an $2m$- dimensional commutative algebra over
${\mathbb R}$, assume that the set  $e_0$, $e_1$, \dots ,
$e_{2m-1}$ is a basis of $A_0$ with the Cayley table:
\[e_ie_j=e_{i\oplus j},\]
where $i\oplus j=i+j$ $\left(mod\ 2m\right)$.

\noindent Algebra $A_0$ has the following matrix representation:
\[e_k\to P_k=P^k_1,\]
where $P_1={\left[p_{ij}\right]}_{2m\times2m}$, $p_{ii+1}=1$ for
$0\le i\le 2m-1$, $p_{2m0}=1$ and $p_{ij}=0$ for the rest of $i,\
j$.

\noindent Let us put

\[{{\mathbf \tau }}^l_0=e_l,\ \ l=0,1,\dots ,2m-1,\]
\[{{\mathbf \tau }}^l_1=e_li \sin  s,\ \ l=0,1,\dots ,2m-1,\]
\[{{\mathbf \tau }}^l_2=e_l\cos s,\ \ l=0,1,\dots ,2m-1,\]
\[{{\mathbf \tau }}^l_{2k}=e_l \cos ks,\ \ {{\mathbf \tau }}^l_{2k+1}=e_li\sin (k+1)s,\ \ l=0,1,\dots ,2m-1,\]
\[k=0,1,2,\dots.\]

\noindent It is easily that ${\mathbf \tau }^0_{2n}{{\mathbf \tau
}}^0_{2k}=\frac{1}{2}\left({{\mathbf \tau
}}^0_{2\left(n-k\right)}+{{\mathbf \tau
}}^0_{2\left(n+k\right)}\right),\ \ n\ge k$,
\[{{\mathbf \tau }}^{l_1}_{2n+1}{{\mathbf \tau }}^{l_2}_{2k+1}=\frac{1}{2}\left({{\mathbf \tau }}^{l_1\bigoplus l_2}_{2\left(n-k\right)}-{{\mathbf \tau }}^{l_1\bigoplus l_2}_{2\left(n+k\right)}\right),\ \ n\ge k, \]
\[{{\mathbf \tau }}^{l_1}_{2n+1}{{\mathbf \tau }}^{l_2}_{2k}=\frac{1}{2}\left({{\mathbf \tau }}^{l_1\bigoplus l_2}_{2\left(n-k\right)+1}+{{\mathbf \tau }}^{l_1\bigoplus l_2}_{2\left(n+k\right)+1}\right),\ \ n\ge k.\]

\noindent Let us introduce the following algebra
\[A=\left\{\sum^{+\infty }_{k=0}{\sum^{2m-1}_{l=0}{\left(a^l_{2k}{{\mathbf \tau }}^l_{2k}+a^l_{2k+1}{{\mathbf \tau }}^l_{2k+1}\right)}}\mathrel{\left|\vphantom{\sum^{+\infty }_{k=0}{\sum^{2m-1}_{l=0}{\left(a^l_{2k}{{\mathbf \tau }}^l_{2k}+a^l_{2k+1}{{\mathbf \tau }}^l_{2k+1}\right)}} a^l_j\in {\mathbb R}}\right.\kern-\nulldelimiterspace}a^l_j\in {\mathbb R}\right\},\]
where $\sum^{+\infty
}_{k=0}{\sum^{2m-1}_{l=0}{\left({\left|a^l_{2k}\right|}^2+{\left|a^l_{2k+1}\right|}^2\right)}}<+\infty
$.

\noindent It is easily verified that $A$ is commutative.

\noindent We consider the subspace $B=\left\{a_0{{\mathbf \tau
}}^1_1+a_1{{\mathbf \tau }}^1_2+a_2{{\mathbf \tau
}}^0_0\mathrel{\left|\vphantom{a_0{{\mathbf \tau
}}^1_1+a_1{{\mathbf \tau }}^1_2+a_2{{\mathbf \tau }}^0_0 {\
a}_i\in {\mathbb R}}\right.\kern-\nulldelimiterspace}{\ a}_i\in
{\mathbb R}\right\}$ of the algebra $A$.

\noindent Let us introduce the function ${\mathbf f}:B\to A$
(${\mathbf f}\left(t,y,z\right)=f\left(e_1\left(t{\cos  s\
}+yi{\sin  s\ }\right)+z\right))$   as follows
\[{\mathbf f}\left(t,y,z\right)=\sum^{+\infty }_{k=0}{\sum^{2m-1}_{l=0}{\left(v^l_{2k}\left(t,y,z\right){{\mathbf \tau }}^l_{2k}+v^l_{2k+1}\left(t,y,z\right){{\mathbf \tau }}^l_{2k+1}\right)}}.\]
The function ${\mathbf f}$\textbf{ }is called $B/A$ differentiable
at ${{\mathbf x}}_0\in B$ if there exists ${{\mathbf f}}^{{\mathbf
'}}\left({{\mathbf x}}_0\right)\in A$ such that for any ${\mathbf
h}\in B$
\[{{\mathbf f}}^{{\mathbf '}}\left({{\mathbf x}}_0\right){\mathbf h}{\mathbf =}{\mathop{\lim }_{\varepsilon \to 0} \frac{{\mathbf f}\left({{\mathbf x}}_0+\varepsilon {\mathbf h}\right)-{\mathbf f}\left({{\mathbf x}}_0\right)}{\varepsilon }\ }\]
In [2] proved that if ${\mathbf f}$ is $B/A$ differentiable, then

\begin{equation}\label{7}
\frac{\partial}{\partial t}{\mathbf f}=e_1{\cos s\
}\frac{\partial}{\partial z}{\mathbf f}
\end{equation}

and

\begin{equation}\label{8}
\frac{\partial}{\partial y}{\mathbf f}=e_1i{\sin  s\
}\frac{\partial}{\partial z}{\mathbf f}.
\end{equation}

In this case all $v^l_{2k}\left(t,y,z\right)$ are solutions of
Eq.(6). Indeed,
\[{\left(\frac{{\partial }^2}{\partial t^2}-\frac{{\partial }^2}{\partial y^2}\right)}^m{\mathbf f}-\frac{{\partial }^{2m}}{\partial z^{2m}}{\mathbf f}=e^{2m}_1{\left({{\cos }^{{\rm 2}} s\ }-{\left(i{\sin  s\ }\right)}^2\right)}^m-1=0.\]
In the sequel we denote by ${\mathbf e}$ the element $e_1$.

\noindent We will seek a solution of Eq.(5) in the following form
\[g_c\left({\mathbf e}\left(t{\cos  s\ }+yi{\sin  s\ }\right)\right)=e^{{\mathbf e}\left(t{\cos  s\ }+yi{\sin  s\ }\right)}\]
Since $f\left({\mathbf e}\left(t{\cos  s\ }+yi{\sin  s\
}\right)+z\right)=g_c({\mathbf e}\left(t{\cos  s\ }+yi{\sin  s\
}\right))e^z$ we have
\[v^l_k\left(t,y,z\right)=u^l_k\left(t,y\right)e^z, l=0,1,\dots ,2m-1, \ k=0,1,2,\dots ,\]
where $g_c\left(t,y\right)=\sum^{+\infty
}_{k=0}{\sum^{2m-1}_{l=0}{\left(u^l_{2k}\left(t,y\right){{\mathbf
\tau }}^l_{2k}+u^l_{2k+1}\left(t,y\right){{\mathbf \tau
}}^l_{2k+1}\right)}}.$

Therefore, we obtain functions $u^l_0\left(t,y\right)$ for $t\ge
\left|y\right|$ from the following equation
\[\sum^{2m-1}_{l=0}{u^l_0\left(t,y\right){{\mathbf \tau }}^l_0}=\sum^{2m-1}_{l=0}{u^l_0\left(t,y\right){{\mathbf e}}^l}\]
\[=\ \frac{1}{2\pi }\int^{\pi }_{-\pi }{e^{{\mathbf e}\left(t{\cos  s\ }+yi{\sin  s\ }\right)}ds}=J_0\left({\mathbf e}i\sqrt{y^2-t^2}\right)=I_0\left({\mathbf e}\sqrt{t^2-y^2}\right),\]
where $I_k$ (resp. $J_k$) is the modified Bessel (resp. Bessel)
function of the first kind and $k$th order [4].

where $I_k$ (resp. $J_k$) is the modified Bessel (resp. Bessel)
function of the first kind and $k$th order [4].

It follows from Eqs.(7),(8) the following Cauchy-Riemann type
conditions

\begin{equation*}
\frac{\partial }{\partial t}u^{l\bigoplus 1}_0=\frac{1}{2}u^l_2,
\end{equation*}

\begin{equation*}
\frac{\partial }{\partial t}u^{l\bigoplus 1}_1=\frac{1}{2}u^l_3,
\end{equation*}

\begin{equation*}
\frac{\partial }{\partial t}u^{l\bigoplus
1}_2=u^l_0+\frac{1}{2}u^l_4,
\end{equation*}

\begin{equation*}
\frac{\partial }{\partial t}u^{l\bigoplus
1}_{2k-1}=\frac{1}{2}\left(u^l_{2k-3}+u^l_{2k+1}\right),
\end{equation*}

\begin{equation}\label{9}
\frac{\partial u^{l\bigoplus 1}_{2k}}{\partial
t}=\frac{1}{2}\left(u^l_{2\left(k-1\right)}+u^l_{2\left(k+1\right)}\right),
\end{equation}
$$k=2,3,\dots ;$$

and

\begin{equation*}
\frac{\partial }{\partial y}u^{l\bigoplus
1}_0=-\frac{1}{2}u^l_1,
\end{equation*}

\begin{equation*}
\frac{\partial u^{l\bigoplus 1}_1}{\partial
y}=u^l_0-\frac{1}{2}u^l_4,
\end{equation*}

\begin{equation*}
\frac{\partial }{\partial y}u^{l\bigoplus 1}_2=-\frac{1}{2}u^l_3,
\end{equation*}

\begin{equation*}
\frac{\partial u^{l\bigoplus 1}_{2k+1}}{\partial
y}=\frac{1}{2}\left(u^l_{2k}-u^l_{2\left(k+2\right)}\right),
\end{equation*}

\begin{equation}\label{10}
\frac{\partial u^{l\bigoplus 1}_{2k+2}}{\partial
y}=\frac{1}{2}\left(u^l_{2k-1}-u^l_{2k+3}\right),
\end{equation}

$$\ k=1,2,\dots .$$

By using Eqs.(9),(10) and functions $u^l_0\left(t,y\right)$, we
can obtain recurrently function $u^l_k\left(t,y\right)$ for any
$k\ge 1,$ which will be used for solution of Eq.(1).

In the sequel, unless otherwise specified, the case where $m=2$ is
studied. In this case $f_s\left(t,y\right)$ is of the following
form $f_s\left(t,x\right)=e^{-t}\left(1+t\right)\delta
\left(t-x\right)$ and hence,
\[g_s\left(t,y\right)=\left(1+t\right)\delta \left(t-y\right).\]
Algebra $A_0$ is as follows
\[A_0{\rm \ }=\left\{a+e_1b+e_2c+e_3d\mathrel{\left|\vphantom{a+e_1b+e_2c+e_3d \ \ a,b,c,d\in {\mathbb R}{\rm \ }}\right.\kern-\nulldelimiterspace}\ \ a,b,c,d\in {\mathbb R}{\rm \ }\right\},\]
where the basis $e_l={{\mathbf e}}^l$, $l=0,1,2,3,$ and ${\mathbf
e}$ has the following matrix representation
\[{\mathbf e}\to P=\left[ \begin{array}{c}
 \begin{array}{cc}
0 & 1 \\
0 & 0 \end{array} \ \ \ \  \begin{array}{cc}
0 & 0 \\
1 & 0 \end{array}
 \\
 \begin{array}{cc}
0 & 0 \\
1 & 0 \end{array} \ \ \ \  \begin{array}{cc}
0 & 1 \\
0 & 0 \end{array}
 \end{array}
\right]{\mathbf .\ }\]
Therefore, we have ${{\mathbf \tau
}}^0_0=1$,  ${{\mathbf \tau }}^0_{2k}={\cos  ks\ }$, ${{\mathbf
\tau }}^l_{2k}={{\mathbf e}}^l{\cos  ks\ }$, ${{\mathbf \tau
}}^l_{2k+1}={{\mathbf e}}^li{\sin  \left(k+1\right)s\ }$,
$l=0,1,2,3$, $\ k=0,1,2,\dots $

Taking into account that $g_c\left({\mathbf e}\left(t{\cos  s\
}+yi{\sin  s\ }\right)\right)=e^{{\mathbf e}\left(t{\cos  s\
}+yi{\sin  s\ }\right)}$, we have
\[u^0_0\left(t,y\right)+{\mathbf e}u^1_0\left(t,y\right)+{{\mathbf e}}^2u^2_0\left(t,y\right)+{{\mathbf e}}^3u^3_0\left(t,y\right)=\ \frac{1}{2\pi }\int^{\pi }_{-\pi }{e^{{\mathbf e}\left(t{\cos  s\ }+yi{\sin  s\ }\right)}ds}\\
=I_0\left({\mathbf e}\sqrt{t^2-y^2}\right).\]
It is easily seen that
\begin{align}
I_0\left({\mathbf e}\sqrt{t^2-y^2}\right)\nonumber \\& = &
\frac{I_0\left(\sqrt{t^2-y^2}\right)
+I_0\left(i\sqrt{t^2-y^2}\right)}{2}+{{\mathbf
e}}^2\left(\frac{I_0\left(\sqrt{t^2-y^2}\right)-I_0\left(i\sqrt{t^2-y^2}\right)}{2}\right)\nonumber\\
& = &
\frac{I_0\left(\sqrt{t^2-y^2}\right)+J_0\left(\sqrt{t^2-y^2}\right)}{2}+{{\mathbf
e}}^2\frac{I_0\left(\sqrt{t^2-y^2}\right)-J_0\left(\sqrt{t^2-y^2}\right)}{2}.\nonumber
\end{align}

Therefore, for $t\ge \left|y\right|$, we have
$u^1_0\left(t,y\right)=u^3_0\left(t,y\right)=0$ and

\begin{align}
u^0_0\left(t,y\right)=\frac{I_0\left(\sqrt{t^2-y^2}\right)+J_0\left(\sqrt{t^2-y^2}\right)}{2},\nonumber\\
u^2_0\left(t,y\right)=\frac{I_0\left(\sqrt{t^2-y^2}\right)-J_0\left(\sqrt{t^2-y^2}\right)}{2}.\nonumber
\end{align}

It follows from the first two equations of (10) that

\begin{align}u^1_1=-2\frac{\partial }{\partial y}u^2_0 = & -\frac{\partial \left[I_0\left(\sqrt{t^2-y^2}\right)-J_0\left(\sqrt{t^2-y^2}\right)\right]}{\partial y} \nonumber \\
= &
\frac{y}{\sqrt{t^2-y^2}}\left(I_1\left(\sqrt{t^2-y^2}\right)+J_1\left(\sqrt{t^2-y^2}\right)\right),\nonumber
\end{align}
\begin{align}u^3_1=-2\frac{\partial }{\partial y}u^0_0= & -\frac{\partial \left[I_0\left(\sqrt{t^2-y^2}\right)+J_0\left(\sqrt{t^2-y^2}\right)\right]}{\partial y} \nonumber \\
= &
\frac{y}{\sqrt{t^2-y^2}}\left(I_1\left(\sqrt{t^2-y^2}\right)-J_1\left(\sqrt{t^2-y^2}\right)\right),\nonumber
\end{align}
\[u^0_1=-2\frac{\partial }{\partial y}u^1_0=0,\]
\[u^2_1=-2\frac{\partial }{\partial y}u^3_0=0.\]

Then it follows from the Cauchy-Riemann type conditions (9) that

\begin{equation*}
u^0_2\left(t,y\right)=2\frac{\partial
u^1_0\left(t,y\right)}{\partial t}=0;
\end{equation*}

\begin{align}
u^1_2\left(t,y\right)= 2\frac{\partial
u^2_0\left(t,y\right)}{\partial t}
=\frac{\partial\left(I_0\left(\sqrt{t^2-y^2}\right)-J_0\left(\sqrt{t^2-y^2}\right)\right)}{\partial
t}\nonumber
\\
=\frac{t}{\sqrt{t^2-y^2}}\left(I_1\left(\sqrt{t^2-y^2}\right)+J_1\left(\sqrt{t^2-y^2}\right)\right);
\nonumber
\end{align}

\[u^2_2\left(t,y\right)=2\frac{\partial u^3_0\left(t,y\right)}{\partial y}=0;\]

\begin{align}
u^3_2\left(t,y\right)=2\frac{\partial
u^0_0\left(t,y\right)}{\partial t}=\frac{\partial
\left[I_0\left(\sqrt{t^2-y^2}\right)+J_0\left(\sqrt{t^2-y^2}\right)\right]}{\partial
t}\nonumber\\
=\frac{t}{\sqrt{t^2-y^2}}\left(I_1\left(\sqrt{t^2-y^2}\right)-J_1\left(\sqrt{t^2-y^2}\right)\right).\nonumber
\end{align}

\noindent Similarly, for $u^l_3$ we have

\begin{align}
u^0_3= & 2\frac{\partial }{\partial t}u^1_1 =2\frac{\partial
}{\partial
t}\left[\frac{y}{\sqrt{t^2-y^2}}\left(I_1\left(\sqrt{t^2-y^2}\right)+J_1\left(\sqrt{t^2-y^2}\right)\right)\right]
\nonumber\\= &
-\frac{2ty}{\sqrt{{\left(t^2-y^2\right)}^3}}\left(I_1\left(\sqrt{t^2-y^2}\right)+J_1\left(\sqrt{t^2-y^2}\right)\right)\nonumber
\\+ &
\frac{ty}{t^2-y^2}\left(I_0\left(\sqrt{t^2-y^2}\right)+I_2\left(\sqrt{t^2-y^2}\right)+J_0\left(\sqrt{t^2-y^2}\right)-J_2\left(\sqrt{t^2-y^2}\right)\right);
\nonumber \end{align}

\begin{align}
u^2_3= & 2\frac{\partial }{\partial t}u^3_1 = 2\frac{\partial
}{\partial
t}\left[\frac{y}{\sqrt{t^2-y^2}}\left(I_1\left(\sqrt{t^2-y^2}\right)-J_1\left(\sqrt{t^2-y^2}\right)\right)\right]
\nonumber \\= &
-\frac{2ty}{\sqrt{{\left(t^2-y^2\right)}^3}}\left(I_1\left(\sqrt{t^2-y^2}\right)-J_1\left(\sqrt{t^2-y^2}\right)\right)\nonumber
\\+ &\frac{2ty}{t^2-y^2}\left(I_0\left(\sqrt{t^2-y^2}\right)+I_2\left(\sqrt{t^2-y^2}\right)-J_0\left(\sqrt{t^2-y^2}\right)+J_2\left(\sqrt{t^2-y^2}\right)\right).\nonumber \end{align}

\noindent It is easily seen that $u^1_3=u^3_3=0.$

\noindent Next, it follows from (9) that
\begin{align}
u^0_4= & 2\frac{\partial u^1_2}{\partial t}-2u^0_0=2\frac{\partial
}{\partial
t}\frac{t}{\sqrt{t^2-y^2}}\left(I_1\left(\sqrt{t^2-y^2}\right)+J_1\left(\sqrt{t^2-y^2}\right)\right)-2u^0_0
\nonumber \\= &
\frac{-2y^2}{\sqrt{{\left(t^2-y^2\right)}^3}}\left(I_1\left(\sqrt{t^2-y^2}\right)+J_1\left(\sqrt{t^2-y^2}\right)\right)\nonumber
\\+ & \frac{t^2}{t^2-y^2}\left(I_0\left(\sqrt{t^2-y^2}\right)+I_2\left(\sqrt{t^2-y^2}\right)+J_0\left(\sqrt{t^2-y^2}\right)-J_2\left(\sqrt{t^2-y^2}\right)\right)\nonumber
\\ & -I_0\left(\sqrt{t^2-y^2}\right)-J_0\left(\sqrt{t^2-y^2}\right);\nonumber \end{align}

\begin{align}
u^2_4= & 2\left({\rm \ }\frac{\partial u^3_2}{\partial
t}-u^2_0\right)=\frac{-2y^2}{\sqrt{{\left(t^2-y^2\right)}^3}}\left(I_1\left(\sqrt{t^2-y^2}\right)-J_1\left(\sqrt{t^2-y^2}\right)\right)\nonumber
\\+ & \frac{t^2}{t^2-y^2}\left(I_0\left(\sqrt{t^2-y^2}\right)+I_2\left(\sqrt{t^2-y^2}\right)-J_0\left(\sqrt{t^2-y^2}\right)+J_2\left(\sqrt{t^2-y^2}\right)\right)\nonumber
\\- & I_0\left(\sqrt{t^2-y^2}\right)+J_0\left(\sqrt{t^2-y^2}\right).\nonumber \end{align}

\noindent Also it is easily verified that $u^1_4=u^3_4=0.$

\noindent By using well-known integrals for Bessel functions
[3-5], we have
\begin{align}
\int^t_{-t}{u^0_0}dy={\sinh  t\ }+{\sin  t\
},\int^t_{-t}{u^2_0}dy={\sinh  t\ }-{\sin  t\ },\
\int^t_{-t}{u^1_1}dy=\int^t_{-t}{u^3_1}dy=0,\nonumber \end{align}
\begin{align}
\int^t_{-t}{u^1_2}dy= & 2\int^t_{-t}{\frac{\partial u^2_0}{\partial t}}dy=2\left(\frac{\partial }{\partial t}\int^t_{-t}{u^2_0}dy-u^2_0\left(t,t\right)-u^2_0(t,-t)\right) \nonumber \\= & 2{\cosh  t\ }-2{\cos  t\ },\nonumber
\end{align}
\begin{align}
\int^t_{-t}{u^3_2}dy= & 2\int^t_{-t}{\frac{\partial u^0_0}{\partial t}}dy=2\left(\frac{\partial }{\partial t}\int^t_{-t}{u^0_0}dy-u^0_0\left(t,t\right)-u^0_0\left(t,-t\right)\right)\nonumber \\= & 2{\cosh  t\ }+2{\cos  t\ }-4.\nonumber
\end{align}

\noindent As an example, we obtain the pdf for the case, where
$m=1$. For this case $e_1=1$ and hence, we can consider functions
$\sum^4_{l=0}{u^l_k\left(t,y\right)}$, $k=0,1,2,\dots $ as
solutions of Eq.(5) for $m=1$.

\noindent For $t\le \left|y\right|$ consider the function
$g\left(t,y\right)=g_c\left(t,y\right)+g_s\left(t,y\right)$ of the
following form:
\begin{align}
g_c\left(t,y\right)= &
\frac{1}{2}\left(u^0_0\left(t,y\right)+u^2_0\left(t,y\right)\right)+\frac{1}{4}\left(u^1_1\left(t,y\right)
+u^3_1\left(t,y\right)+u^1_2\left(t,y\right)+u^3_2\left(t,y\right)\right)\nonumber
\\=
&\frac{I_0\left(\sqrt{t^2-y^2}\right)}{2}+\frac{t+y}{2\sqrt{t^2-y^2}}I_1\left(\sqrt{t^2-y^2}\right)\nonumber
\end{align} and $g_s\left(t,y\right)=\delta \left(t-y\right)$.

\noindent It is easily seen that function $g_c\left(t,y\right)$ is
a solution of the equation for $t<y$

\begin{align}
\left(\frac{{\partial }^2}{\partial t^2}-\frac{{\partial
}^2}{\partial y^2}\right)g\left(t,y\right)-g\left(t,y\right)=0,
\end{align}

In addition, we have ${\mathop{\lim }_{y\uparrow t}
g_c\left(t,y\right)\ }=\frac{1}{2}\left(1+t\right)$ and
${\mathop{\lim }_{y\downarrow -t} g_c\left(t,y\right)\
}=\frac{1}{2}.$

\noindent To avoid cumbersome calculations we put $v=1$.

\noindent Therefore, $f\left(t,x\right)=e^{-t}g\left(t,x\right)$
is a solution of the equation:

\begin{align}
\left(\frac{\partial }{\partial t}-\frac{\partial }{\partial
x}+1\right)\left(\frac{\partial }{\partial t}+\frac{\partial
}{\partial x}+1\right)f_c\left(t,x\right)-f_c\left(t,x\right)=0,
\end{align}

\[f_s\left(t,x\right)=\delta \left(t-x\right)e^{-t}.\]

\noindent In addition, $f_c\left(t,x\right)$ satisfies the
following conditions:
\[{\mathop{\lim }_{x\uparrow t} f_c\left(t,x\right)\ }=\frac{1}{2}\left(e^{-t}+te^{-t}\right),\ \ {\mathop{\lim }_{x\downarrow -t} f_c\left(t,x\right)\ }=\frac{1}{2}e^{-t},\]
and for all $t>0$ we have
$\int^t_{-t}{f\left(t,x\right)}dx=1$.\textit{}

\noindent For a small $\varepsilon >0$ consider the probability
$P\left\{0<t-x\left(t\right)<\varepsilon \right\}$.

\noindent Let us verify that ${\mathop{\lim }_{x\uparrow t}
f_c\left(t,x\right)\ }=\mathop{{\rm lim}}_{\varepsilon \downarrow
0}\ \frac{P\left\{0<t-x\left(t\right)<\varepsilon
\right\}}{\varepsilon }$, i.e.,
$$\mathop{{\rm lim}}_{\varepsilon
\downarrow 0}\ \frac{P\left\{0<t-x\left(t\right)<\varepsilon
\right\}}{\varepsilon }=\frac{1}{2}\left(e^{-t}+te^{-t}\right).$$

\noindent Indeed, it is easily seen that

\begin{align}P\left\{0<t-x\left(t\right)\le \varepsilon \right\}=P\left\{t-\frac{\varepsilon }{2}\le {\theta }_1<t\right\}+\int^t_0{P\left\{{\theta }_3\ge t-u,{\theta }_2\le \frac{\varepsilon }{2},{\theta }_1\in du\right\}}\nonumber\\
+o(\varepsilon ),\nonumber \end{align} where ${\theta }_i$,
$i=1,2,3$ are independent exponentially distributed random
variables.

\noindent The random variable ${\theta }_1$ is time of the first
velocity alternation, ${\theta }_2$ is time between the first and
the second velocity alternations and ${\theta }_3$ is time between
the second and the third velocity alternations.

\noindent We have that $P\left\{t-\frac{\varepsilon }{2}\le
{\theta }_1<t\right\}=e^{-t+\frac{\varepsilon }{2}}-e^{-t}$ and as
it easy to calculate
\begin{align}\int^t_0{P\left\{{\theta }_3\ge t-u,\ {\theta }_2\le \frac{\varepsilon }{2},{\theta }_1\in du\right\}} = & \left(1-e^{-\frac{\varepsilon }{2}}\right)\int^t_0{e^{-t+u}e^{-u}du}\nonumber \\
=& \left(1-e^{-\frac{\varepsilon }{2}}\right)te^{-t}.\nonumber
\end{align}
Whence, it is easily verified that $\mathop{{\rm
lim}}_{\varepsilon \downarrow 0}\
\frac{P\left\{0<t-x\left(t\right)<\varepsilon
\right\}}{\varepsilon }=\frac{1}{2}\left(e^{-t}+te^{-t}\right)$.

\noindent Similarly, $P\left\{t+x\left(t\right)\le \varepsilon
\right\}=P\left\{{\theta }_2\ge t-\frac{\varepsilon }{2},{\theta
}_1\le \frac{\varepsilon }{2}\right\}+o\left(\varepsilon \right).$
This implies that
\[\mathop{{\rm lim}}_{\varepsilon \downarrow 0}\ \frac{P\left\{t+x\left(t\right)<\varepsilon \right\}}{\varepsilon }=\frac{1}{2}e^{-t}={\mathop{\lim }_{x\downarrow -t} f_c\left(t,x\right)\ }.\]

Therefore, $f_c\left(t,x\right)$ is a solution of the Goursat
problem for the linear second order hyperbolic equation that
ensures the uniqueness of the solution of Eq.(12) with conditions
(4). It means that $f\left(t,x\right)$ is the pdf of the
particle's position for $m=1$.

\noindent It is relevant to remark that function
$f\left(t,x\right)$ coincides with the result obtained in [5].

Now, we turn to the case $m=2$ and continue to calculate integrals
of $u^l_k$.

\noindent It follows from $u^0_4=2\frac{\partial u^1_2}{\partial
t}-2u^0_0$ that
\begin{align}\int^t_{-t}{u^0_4}dy = & 2\left(\frac{\partial }{\partial t}\int^t_{-t}{u^1_2}dy-u^1_2\left(t,t\right)-u^1_2\left(t,-t\right)\right)-2{\sinh  t\ }-2{\sin  t\ }
\nonumber \\ = & 4\left({\sinh  t\ }+{\sin  t\ }-t\right)-2{\sinh
t\ }-2{\sin t\ }=2{\sinh  t\ }+2{\sin  t\
}-4t.\nonumber\end{align}

Next, it follows from $u^2_4=2\frac{\partial u^3_2}{\partial
t}-2u^2_0$ that
\begin{align}\int^t_{-t}{u^2_4}dy= & 2\left(\frac{\partial }{\partial t}\int^t_{-t}{u^3_2}dy-u^3_2\left(t,t\right)-u^3_2\left(t,-t\right)\right)-2{\sinh  t\ }+2{\sin  t\ }\nonumber \\
= & 4{\sinh  t\ }-4{\sin  t\ }-2{\sinh  t\ }+2{\sin  t\ }\
=2{\sinh  t\ }-2{\sin  t\ }.\nonumber
\end{align}

\noindent For $t\le \left|y\right|$ we introduce the function
$g\left(t,y\right)=g_c\left(t,y\right)+g_s\left(t,y\right)$, where

\begin{equation*}
g_c\left(t,y\right)=\frac{1}{2}u^2_0\left(t,y\right)+\frac{1}{4}\left(u^1_1\left(t,y\right)+u^3_1\left(t,y\right)+u^1_2\left(t,y\right)+u^3_2\left(t,y\right)+u^0_4\left(t,y\right)\right),
\end{equation*}
\begin{equation}g_s\left(t,y\right)=\delta \left(t-y\right)+t\delta
\left(t-y\right).
\end{equation}

\noindent By construction, the function $g_c\left(t,y\right)$ is a
solution of the following equation

\begin{equation}
{\left(\frac{{\partial }^2}{\partial t^2}-\frac{{\partial
}^2}{\partial y^2}\right)}^2g\left(t,y\right)-g\left(t,y\right)=0.
\end{equation}

Therefore, the function
$f_c\left(t,x\right)=e^{-t}g_c\left(t,x\right)$ is a solution of
Eq.(1) for $m=2$ ($\lambda =v=1$).

We put
$f\left(t,x\right)=f_c\left(t,x\right)+e^{-t}g_s\left(t,x\right)$.
Taking into account the values of integrals of functions, which
are involved in the expression for $g_c\left(t,y\right)$, we have
that $\int^t_{-t}{f\left(t,x\right)}dx=1$, for all $t\ge 0.$

Let us prove that ${\mathop{\lim }_{x\uparrow t}
f_c\left(t,x\right)\ }=\mathop{{\rm lim}}_{\varepsilon \downarrow
0}\ \frac{P\left\{0<t-x\left(t\right)<\varepsilon
\right\}}{\varepsilon }$ and ${\mathop{\lim }_{x\downarrow -t}
f_c\left(t,x\right)\ }=\mathop{{\rm lim}}_{\varepsilon \downarrow
0}\ \frac{P\left\{t+x\left(t\right)<\varepsilon
\right\}}{\varepsilon }$.

\noindent It follows from Lemma 2 that for $m=2$ we have

\begin{equation*}
\mathop{{\rm lim}}_{\varepsilon \downarrow 0}\
\frac{P\left\{0<t-x\left(t\right)<\varepsilon
\right\}}{\varepsilon }=\frac{1}{2}te^{-t}.
\end{equation*}

and
\[\mathop{{\rm lim}}_{\varepsilon \downarrow 0}\ \frac{P\left\{t+x\left(t\right)<\varepsilon \right\}}{\varepsilon }=0.\]
It easily verified that $\mathop{{\rm lim}}_{y\uparrow t}{\
u}^0_4\left(t,y\right)=0$, $\mathop{{\rm lim}}_{y\uparrow t}{\
u}^2_0\left(t,y\right)=0$ and consequently

\begin{align}
\mathop{{\rm lim}}_{y\uparrow t}\
g_c\left(t,y\right)=\lim_{y\uparrow
t}\frac{t+y}{2\sqrt{t^2-y^2}}I_1\left(\sqrt{t^2-y^2}\right)=\frac{t}{2},
\nonumber
\\ \lim_{y\downarrow -t}g_c\left(t,y\right)=\mathop{{\rm
lim}}_{y\downarrow
-t}\frac{t+y}{2\sqrt{t^2-y^2}}I_1\left(\sqrt{t^2-y^2}\right)=0.
\end{align}

Thus,

\begin{align}
\mathop{{\rm lim}}_{x\uparrow t}\
f_c\left(t,x\right)=\frac{1}{2}te^{-t}=\mathop{{\rm
lim}}_{\varepsilon \downarrow 0}\
\frac{P\left\{t-x\left(t\right)<\varepsilon \right\}}{\varepsilon
}\ ,\nonumber \\
{\mathop{\lim }_{x\downarrow -t}
f_c\left(t,x\right)\ }=0=\mathop{{\rm lim}}_{\varepsilon
\downarrow 0}\ \frac{P\left\{t+x\left(t\right)<\varepsilon
\right\}}{\varepsilon }.
\end{align}

\noindent Let us show that conditions (15) with the condition
$\int^t_{-t}{g\left(t,y\right)e^{-t}}dx=1$ insure the uniqueness
of the solution $g_c\left(t,y\right)$ for Eq.(14) and
consequently, the uniqueness of the solution $f_c\left(t,x\right)$
of Eq.(12).

\noindent It is easily seen that each solution of Eq.(11) is a
solution of Eq.(14). By changing the variables $s=t+y$, $p=t-y$,
we reduce Eq.(14) to

\begin{equation}
\frac{{\partial }^4}{\partial s^2\partial
p^2}G\left(s,p\right)-G\left(s,p\right)=0.
\end{equation}

Passing to the Fourier transform $\hat{G}\left(s,\alpha
\right)=\int^{\infty }_0{G\left(s,p\right)e^{i\alpha p}dp}$ in
Eq.(17), we get the ordinary differential equation of order 4.
Taking into account that $\mathop{{\rm lim}}_{y\downarrow
-t}g_c\left(t,y\right)=0$, we have

\begin{equation}
\hat{G}\left(0,\alpha \right)=0.
\end{equation}

Hence, at most four independent solutions of the ordinary
differential equation satisfy the initial condition (18) for each
$\alpha $. Passing to the inverse Fourier transform, we have four
independent solutions of Eq.(14) with the condition ${\mathop{\lim
}_{x\downarrow -t} g_c\left(t,x\right)\ }=0$ and just two of them
satisfy Eq.(14) but not Eq.(11). By construction, one of these
solutions $g_c\left(t,y\right)$ is given by Eq.(13). As another
solution we can take
\[g_2(t,y)={\ u}^2_0\left(t,y\right)+{\ u}^0_4\left(t,y\right).\]
It is easily verified that no linear combination $c(t,y)$ of
functions $g_c\left(t,y\right)$ and $g_2(t,x)$ satisfies
conditions (16) and
$\int^t_{-t}{\left(c\left(t,x\right)+g_s\left(t,y\right)\right)e^{-t}}dx=1$
for all $t>0$, but solution $g_c\left(t,y\right)$.

\noindent Therefore, the function $f\left(t,x\right)$ is the pdf
of the particle position at time $t$ for $m=2$, $v=\lambda =1$ and
has the following form
\begin{align} f\left(t,x\right)= & -\frac{J_0\left(\sqrt{t^2-x^2}\right)}{2}e^{-t}+\frac{\left(t+x\right)e^{-t}}{2\sqrt{t^2-x^2}}I_1\left(\sqrt{t^2-x^2}\right)\nonumber\\
-& \frac{x^2e^{-t}}{2\sqrt{{\left(t^2-x^2\right)}^3}}\left(I_1\left(\sqrt{t^2-x^2}\right)+J_1\left(\sqrt{t^2-x^2}\right)\right)\nonumber\\
+&
\frac{t^2e^{-t}}{4\left(t^2-x^2\right)}\left(I_0\left(\sqrt{t^2-x^2}\right)+I_2\left(\sqrt{t^2-x^2}\right)
+J_0\left(\sqrt{t^2-x^2}\right)-J_2\left(\sqrt{t^2-x^2}\right)\right)\nonumber\\
+&
\delta\left(t-x\right)e^{-t}+t\delta\left(t-x\right)e^{-t}.\nonumber\end{align}
In much the same way as the pdf $f\left(t,x\right)$ of the
particle position for $m=2$ was obtained we can also get solutions
of Eq.(1) with conditions (2) and (4) for each $m>2$.

\noindent

\noindent \eject

\end{document}